\documentclass[12pt]{amsart}
\usepackage{amssymb}
\newcommand{\tensor}{\otimes}

\renewcommand{\O}{\mathcal O}
\newcommand{\A}{\mathbb A}
\newcommand{\R}{\mathbb R}
\newcommand{\RR}{\mathbb R}

\newcommand{\Qbar}{\hat Q}
\newcommand{\N}{\mathbb N}
\newcommand{\NN}{\mathbb N}
\newcommand{\Z}{\mathbb Z}
\newcommand{\C}{\mathbb C}
\newcommand{\p}{\bf p}%
\newcommand{\sdr}{\ltimes}

\newcommand{\Aut}{\mathop{\rm Aut}}

\newcommand{\codim}{\mathop{\rm codim}}
\newcommand{\Diff}{\mathop{\rm Diff}}

\newcommand{\dd}[1]{\frac{\partial}{\partial #1}}

\newcommand{\mk}{{N(m,n,k)}}
\newcommand{\mkn}{{[\binom{m+k}{m}-(m+1)]\cdot(n-m)}}
\newcommand{\B}{{\rm B}}

\newcommand{\dist}{\rm dist}
\newcommand{\id}{\rm id}
\newcommand{\Auth}{{\rm Aut}_{hol}}

\newcommand{\ha}{\hookrightarrow}
\newcommand{\f}{\varphi}
\theoremstyle{plain}%

\newtheorem{corollary}{Corollary}
\newtheorem{lemma}{Lemma}
\newtheorem{proposition}{Proposition}
\newtheorem*{sublemma}{Sublemma}
\newtheorem{theorem}{Theorem}
\newtheorem{theoremI}{Theorem}

\theoremstyle{definition}

\newtheorem{definition}{Definition}
\newtheorem*{remark}{Remark}
\newtheorem{example}{Example}

\newtheorem{question}{Question}
\begin{document}
\title{Subvarieties of $\C^n$ with Non-Extendible Automorphisms}
\footnote{All authors were partially supported by Schweizerische
Nationalfonds (SNF), the first author was also partially supported
by the ``Freiwillige Akademische Gesellschaft''.}
\author{Harm Derksen}
\address{Derksen, Northeastern University\\ hderksen@math.neu.edu}
\author{Frank Kutzschebauch}
\address{Kutzschebauch, Universit\"at Basel\\ kutzsche@math.unibas.ch}
\author{J\"org Winkelmann}
\address{Winkelmann, Universit\"at Basel\\ jwinkel@member.ams.org}
%
\subjclass{14E09, 14L30, 32H02, 32M05}
\begin{abstract}
We investigate algebraic and analytic subvarieties of $\C^n$
with automorphisms which can not be extended to the
ambient space.

\end{abstract}
\maketitle
\section{Introduction}
If $Z$ is an affine variety
and $n$ is a sufficiently large number, then any two embeddings
$i,j:Z\to\A^n$ are equivalent in the sense that there exists an
automorphism $\phi$ of $\A^n$ such that $\phi\circ i= j$
(see \cite J \cite K \cite S). In particular, in this case
every automorphism of $Z$ extends to the whole of $\A^n$.
In the case where $Z$ is smooth it suffices to take $n>2\dim Z+1$.

Thus for an algebraic subvariety $Z\subset\A^n$ of high codimension every
automorphism extends to an automorphism of $\A^n$.

This raises three questions:
\begin{enumerate}
\item Is a similar statement true for subvarieties of low codimension,
e.g.~hypersurfaces?
\item
Does a similar statement hold in the holomorphic category?
\item
Assume that $Z\subset\A^n$ is a subvariety such that every single
automorphism of $Z$ extends to $\A^n$. Does this imply that there
is an extension of the action of the group $\Aut(Z)$
to an action on $\A^n$?
\end{enumerate}

The purpose of this article is to provide negative answers to all
three questions.

First, we prove that there exist
irreducible algebraic hypersurfaces admitting a $(\C,+)$-action
which does not extend.
\begin{theoremI}
There exists a smooth irreducible algebraic
hypersurface $H\subset \A^5$ and an algebraic action $\mu$ of the additive
group $(\C,+)$ on $H$ such that for all $t\ne 0$
there exists neither an algebraic nor a holomorphic
automorphism $\phi$ of $\A^5$ with $\phi|_H=\mu(t)$.
\end{theoremI}

The key idea for our construction
is to choose a hypersurface in $\A^n$ such
that the complement has a small automorphism group and then
to construct an automorphism on the hypersurface based on something which
does not extend from the hypersurface to the whole space.
For this it is crucial that there are
free actions of algebraic groups $G$ on the affine space $\A^n$ for which
there exist $G$-invariant subvarieties $W$ such that the restriction
map $\O(\A^n)^G\to\O(W)^G$ is not surjective (\cite{JW1}).

It should be noted that this strategy can not work  for subvarieties
of higher codimension on $\A^n$, because for a subvariety $Z$ in
$\A^n$ of codimension at least two the automorphism group of the
complement is always quite large. In fact the group of all
automorphisms of $\A^n$ fixing every point in $Z$ acts
\begin{em} transitively \end{em} on the complement $\A^n\setminus Z$
as soon as $\codim(Z)\ge 2$ (see  \cite{JW2}).

In the second part, we prove that there does not exist
any effective differentiable, holomorphic or algebraic
action of the group  $\Aut(\C^*\times\C^*)$ on the affine space
$\C^n$.
\begin{theoremI}\label{harmful-theo}
Let $k$ be an algebraically closed field,
$K$ a field, $n\in\N$
and let $G$ denote the group given by the
semidirect product $SL_2(\Z)\sdr_\rho (k^*\times k^*)$
with
\[
\rho\begin{pmatrix} a & b \\ c & d \\ \end{pmatrix}(z,w)=
(z^aw^b, z^cw^d)
\]

Then there does not exist any injective group homomorphism from
$G$  either into the group $\Diff(\RR^n)$ of diffeomorphisms of $\RR^n$
or into the group $\Aut_K(\A^n)$ of $K$-automorphisms of the affine
space $\A^n$.
\end{theoremI}

Hence, whenever $i:Z=\C^*\times\C^*\hookrightarrow\C^n$ has the
property that every automorphism of $Z$ extends to $\C^n$,
there must be a non-trivial non-split short exact sequence
of groups
\[
1 \to L \to \tilde G \stackrel{\alpha}{\to} \Aut(Z) \to 1
\]
such that for every $g\in\tilde G$ the restriction $\mu(g)|_Z$
coincides with the canonical action of $\alpha(g)\in \Aut(Z)$ on
$i(Z)$.

A similar result holds in positive characteristic.

Finally, in the third section,
we discuss the situation in the  analytic category.
Unlike the algebraic situation there are analytic subvarieties
of high codimension such that no non-trivial automorphism extends.

The first such examples are due to Rosay and Rudin  \cite{RR1}
who proved that
for $n\ge 2$ there exist \begin{em} discrete \end{em} subsets $S$ of
$\C^n$ such that no non-trivial permutation of $S$ extends to
an automorphism of $\C^n$.

Buzzard and Forn\ae ss \cite{BF} proved that given a hypersurface
$X$ in $\C^N$
there exists an embedding $j:X\to\C^N$ such that the complement
$\C^N\setminus X$ is
hyperbolic. This implies that at most countably many automorphisms of
$X$ can be extended to automorphisms of $\C^N$.

Thus the existence of analytic subvarieties with non-extendible
automorphisms is well-known for proper analytic subvarieties of
maximal or minimal codimension.

Our contribution is to provide a result about
analytic subvarieties of arbitrary
codimension.

We prove the following:
\begin{theoremI}
For every non-finite
analytic subvariety $X\subset\C^n$, every Lie group $G$ and
every effective $G$-action on $X$ there exists an embedding
$j:X\to\C^n$ such that for no element $g\in G\setminus\{e\}$ can the induced
automorphism of $j(X)$ be extended to an automorphism of $\C^n$.
\end{theoremI}

As a consequence we obtain the following result:
\begin{corollary}
Let $X$ be an infinite Stein manifold such that $\Aut(X)$ is a Lie
group.
Then there exists an embedding $j$ of $X$ into some $\C^N$ such that
no non-trivial automorphism of $X$ extends to $\C^N$.
\end{corollary}
The condition that $\Aut(X)$ is a Lie group holds in particular
if $X$ is hyperbolic, e.g.~if the universal covering of $X$ is
biholomorphic to a bounded domain.

\section{Hypersurfaces}
In this section we investigate the algebraic subvarieties of
low codimension, concentrating on hypersurfaces.
Varieties, functions, maps, group actions etc.
are assumed to be algebraic over some algebraically closed ground
field $k$.
Unless explicitly stated otherwise,
this ground field $k$ may have positive characteristic.
\subsection{Basic tools}
We start with  some basic observations
which will be used in later constructions.
\begin{lemma}\label{ob-equiv}
Let $\pi:X\to Y$ be a separable surjective morphism between
irreducible normal algebraic varieties defined over some ground field $k$.
Let $\phi\in \Aut_k(X)$ and assume that $\pi\circ\phi|_F$ is constant for
every $\pi$-fiber $F$.

Then $\phi$ induces a $k$-automorphism $\phi'$ of $Y$ such that
$\pi\circ\phi=\phi'\circ \pi$.
\end{lemma}
\begin{proof}
Since $\pi$ is separable and dominant, the function field $k(Y)$ is
isomorphic to the subfield $L$ of $k(X)$ containing all those
rational functions on $X$ which are constant along the $\pi$-fibers.
Hence $\phi$ induces an automorphism of the function field $k(Y)$.
On the other hand, $\phi$ induces a permutation of the points of $Y$,
because $\pi$ is surjective, and $\pi\circ\phi$ is constant
along the $\pi$-fibers.
Since $Y$ is normal, a rational function $f$ on $Y$ is regular in a
given point $y$ if and only if $f$ has no pole in $y$.
Using this fact, it follows that $\phi$ induces a regular automorphism
of $Y$.
\end{proof}

\begin{lemma}\label{ob-fct}
Let $H$ be an irreducible hypersurface in the affine space $\A^n$
defined by a (regular) function $f$ and $\Omega=\A^n\setminus H$.

Then the map $f|_{\Omega}:\Omega\to\A^1\setminus\{0\}$ is equivariant
for the whole (algebraic) automorphism group $\Aut(\Omega)$.
\end{lemma}
\begin{proof}
Let $g\in k[\Omega]^*$ and let $n$ denote the multiplicity of $g$
along $H$. Then $gf^{-n}\in k[\A^n]^*=k^*$, hence $g=\alpha f^n$ for
some constant $\alpha\in k^*$. It follows that two points
$p,q\in\Omega$ are in the same $f$-fiber if and only if $g(p)=g(q)$
for all $g\in k[\Omega]^*$. Thus the equivalence relation on $\Omega$
defined by $f$ is natural and must be preserved by all automorphisms
of $\Omega$.
Thus it follows from lemma~\ref{ob-equiv} that every
$\phi\in\Aut(\Omega)$ induces an automorphism $f_*\phi$ of
$\A^1\setminus\{0\}$ such that
$f\circ\phi=(f_*\phi)\circ f$.
\end{proof}

\begin{lemma}\label{ob-dcross}
Let $P$ be a polynomial automorphism of $\A^2$ stabilizing the set
$C=\{(x,y):xy=0\}$. Then either $P(x,y)=(\alpha x,\beta y)$ or
$P(x,y)=(\alpha y,\beta x)$ (with $\alpha,\beta\in k^*$).
\end{lemma}
\begin{proof}
Let $\tau(x,y)=(y,x)$. Then either $P$ or $P\circ\tau$ stabilize both
irreducible components of $C$. By the preceding lemma this
implies that $P$  resp.~$P\circ\tau$ is simultaneously equivariant for
both projections $(x,y)\mapsto x$ and $(x,y)\mapsto y$. From this fact
the statement is easily deduced.
\end{proof}
\subsection{Basic examples}
Here we collect some basic examples which are not irreducible.
\begin{example}
Let $k$ be a field, $z\in k\setminus\{0,1,-1\}$ and $S=\{0,1,z\}$.
Then $\sigma(0)=0$, $\sigma(1)=z$, $\sigma(z)=1$ defines
a permutation of $S$ which does not extend to an automorphism
of $\A^1$.
\end{example}
\begin{proof}
This is immediate, since every automorphism of $\A^1$ is
affine-linear.
\end{proof}
\begin{example}
Let $C\subset\A^2$ be the reducible curve defined by
$C=\{(x,y):x(xy-1)=0\}$. The
action of the additive group $G_a$ given by
$\mu_t(x,y)=(x,y)$ for $xy=1$  and $\mu_t%
(0,y) = (0,y+t)$ can not be extended to an algebraic $G_a$-action
on $\A^2$ although there does exist a holomorphic extension
in the case $k=\C$.
\end{example}
\begin{proof}
Note that the two irreducible components of $C$ are
non-isomorphic.
Hence every automorphism of $\A^2$ stabilizing $C$ must stabilize both
irreducible components separately. It follows that  such an algebraic
automorphism must be equivariant for both $(x,y)\mapsto x$ and
$(x,y)\mapsto xy$. From this it is easily deduced that
any such automorphism is of the form
$(x,y)\mapsto (\alpha x,\frac{1}{\alpha}y)$ for some $\alpha\in k^*$.

But an automorphism of this form can not restrict to
$(0,y)\mapsto (0,y+t)$ on $\{(x,y):x=0\}$.

On the other hand, a holomorphic extension is given by
\[
t:(x,y)\mapsto (x,y e^{-tx} - \frac{e^{-tx}-1}{x}).
\]
\end{proof}
\begin{example}
Let $C\subset \A^2$ be the reducible curve defined by
$C=\{(x,y):xy(xy-1)=0\}$. Let $\phi_0$ denote the algebraic
automorphism of $C$ given by $\phi_0(x,y)=(y,x)$ for $xy=1$ and
$\phi_0(x,y)=(x,y)$ for $xy=0$.

For $k=\C$ there does not exist a homeomorphism $\phi$ of
$\A^2(k)=\C^2$ with
$\phi|_C=\phi_0$.
\end{example}
\begin{proof}
Let $D=\{(x,y):xy=0\}$ and $\Omega=\C^2\setminus D$. Consider
$\zeta_\epsilon:S^1\to\Omega$ given by
$\zeta_\epsilon(z)=(z,\epsilon)$ (with $\epsilon\ne 0$). If $\phi$ is
an homeomorphism of $\C^2$ with $\phi|_D=id$, then $\lim_{\epsilon\to
0}\zeta_\epsilon=(id,0)$ implies that $\phi$ stabilizes the element of
$\pi_1(\Omega)$ corresponding to $\zeta_\epsilon$. Now
$\pi_1(\Omega)\simeq\Z^2$
is generated by this element and similar curves around the
$\{y=0\}$-line. This implies that such a homeomorphism $\phi$ must
induce the identity map on $\pi_1(\Omega)$. Since $\pi_1(Q)$ embeds
into $\pi_1(\Omega)$ for $Q=\{(x,y):xy=1\}$, it follows that $\phi|_Q$
must induce the identity map on $\pi_1(Q)$. Therefore $\phi|_Q$ can
not equal the map $(x,y)\mapsto (y,x)$ on $Q$ for any homeomorphism
$\phi:\C^2\to\C^2$ with $\phi|_D=id_D$.
\end{proof}

\begin{example}
Let $(n,m)\ne(1,1)$ be a pair of coprime natural numbers,
$C=\{(x,y):x^ny^m=1\}$ and $\phi\in \Aut(C)$ be given by
$\phi(z,w)=(1/z,1/w)$.

Then $C$ is a smooth connected curve and $\phi$ is an automorphism of
$C$ which can not be extended to an automorphism of $\A^2$.
\end{example}
\begin{proof}
Due to lemma~\ref{ob-fct} the map $(x,y)\mapsto x^ny^m-1$ is
equivariant for every automorphism $\phi\in\Aut(\A^2)$
which stabilizes $C$.
In particular such an automorphism stabilizes the set
$\{x^ny^m=0\}=\{xy=0\}$, since this is the only reducible
fiber of $(x,y)\mapsto x^ny^m-1$. Hence
lemma \ref{ob-dcross}
 implies that any automorphism of $\A^2$ stabilizing $C$
must be described as either $(\alpha x,\beta y)$ or $(\alpha y,\beta x)$. The
second type doesn't stabilize $C$ and the former doesn't restrict to
$\phi$.
\end{proof}
\begin{example}
Let $H=\{(x,y,z):xy=1\}$ and $\phi:(x,y,z)\mapsto (x,y,xz)$.
Then $\phi|_H$ is an automorphism of $H$ which can not be extended to
an automorphism of $\A^3$.
\end{example}
\begin{proof}
Assume that $\phi|_H$ extends to a polynomial automorphism $P$ of $\A^3$.
Let $\Omega=\{(x,y,z):xy\ne 1\}$ and $U=\{(x,y):xy\ne 1\}$ and
consider the morphisms
\[
\Omega\buildrel \pi\over\longrightarrow U
\buildrel\rho\over\longrightarrow \A^*
\]
given by $\pi(x,y,z)=(x,y)$ and $\rho(x,y)=xy-1$.
By lemma \ref{ob-fct} both $\rho$ and $\rho\circ\pi$ are equivariant for
all automorphisms of $U$ resp.\ $\Omega$.
Now all the $\pi$-fibers
 are lines and the generic $\rho$-fiber is $\A^*$.
Since every morphism from $\A$ to $\A^*$ is constant, it follows
(with the help of lemma~\ref{ob-equiv})
that $\pi$ is equivariant as well.

The map $\rho$ being equivariant implies in particular that
the special fiber $\rho^{-1}(-1)$ is invariant under all automorphisms
of $\A^2$ stabilizing $U$.
It follows that any automorphism $Q$ of $\A^2$ stabilizing $U$ must also
stabilize $\{(x,y)\in\A^2:xy=0\}$.
Lemma~\ref{ob-dcross} thus implies that there exist $\alpha,\beta\in
K$
such that either $Q(x,y)=(\alpha x,\beta y)$ or
$Q(x,y)=(\beta y,\alpha x)$.

However, $Q$ is supposed to fix $\{xy=1\}$ pointwise.
This forces $\alpha=\beta=1$ and $Q=\id_{\A^2}$.

Thus an algebraic automorphism $P$ of $\A^3$ extending $\phi|_H$ can
only act along the $\pi$-fibers and therefore can be written in the
form
\[
P:(x,y,z)\mapsto (x,y,g(x,y,z))
\]
for some function $g$. The determinant of the Jacobian of such a map equals
$\frac{\partial g}{\partial z}$. As a nowhere vanishing regular
function it must be constant. This implies that $g$ can be
written in the form
$g(x,y,z)=g_0(x,y,z^p)+\alpha z$ with $\alpha\in k^*$,
$p=char(k)$ and $g_0\in k[X_1,X_2,X_3]$. But now
$g|_{\{xy=1\}}$ can not equal $xz$ and hence we obtained a contradiction to
the assumption that $\phi_H$ extends to an algebraic automorphism $P$
of $\A^3$.
\end{proof}

\subsection{Main hypersurface example}

\begin{theorem}
There exists a smooth irreducible algebraic
hypersurface $H\subset \A^5$ and an algebraic action $\mu$ of the additive
group $(\C,+)$ on $H$ such that for all $t\ne 0$
there exists neither an algebraic nor a holomorphic
automorphism $\phi$ of $\A^5$ with $\phi|_H=\mu(t)$.
\end{theorem}

\begin{proof}
In \cite{JW1} it is shown that
there exists an algebraic $\C$-principal bundle
$\pi:\C^5\to X$  with $X=Q\setminus (S\cup E)$, where $Q$ is a smooth
projective
quadric, $S$ is a smooth hypersurface and $E$ a
two-dimensional smooth subvariety of $Q$
which  intersects $S$ transversally.
Let $Q_1$ denote the variety obtained by $Q$ by blowing-up $E$.
We may, if necessary, blow-up
$Q_1$ again and thereby assume that there is a projective manifold
$\Qbar$, a divisor $D$ with simple normal crossings as its only
singularities, an irreducible component $D_0$ of $D$ and a
birational connected surjective morphism $\tau:\Qbar\to Q$
with $\tau(D_0)=E$ inducing an isomorphism $\Qbar\setminus D\simeq
Q\setminus (S\cup E)$. Let $F$ be a generic fiber of $\tau|_{D_0}\to
E$
(i.e. $F$ is smooth and connected).
Now we fix a very ample line bundle $\O(1)$ on $\Qbar$ and choose
$n\in\N$ such that $\O(n)\tensor K_{\Qbar}$ is ample. By Bertini's
theorem, for every smooth submanifold $Z\subset \Qbar$ there is a dense
open subset in the linear system $|\O(n)|$ such that every
divisor therein intersects $Z$ transversally.
Let $(D_i)_{i\in I}$ be the family of irreducible components of $D$.
Applying Bertini's theorem to all $\cap_{j\in J}D_j$ with $J\subset
I$, $F$ and $\Qbar$ itself we may conclude that there is a very ample
divisor $H$ on $\Qbar$ such that $D\cup H$ is again a divisor with
only simple normal crossings as singularities and furthermore such
that $H$ intersects $F$ transversally.
Let $H_0=\tau(H)\subset Q$. Since $H$ intersects $F$ transversally, it
is clear that $E\subset H_0$. Now $H_0\setminus S$ is affine and
$E\setminus S$ is a hypersurface in $H_0$. It follows that there
exist rational function on $H_0$ which are regular on
$H_0\setminus (E\cup S)$ and have poles of arbitrarily high
multiplicity on $E\setminus S$.
Fix such a function $f$ of a sufficiently high pole order.
We will see later, at the end of the proof, how high this multiplicity has to
be.
Now let $\rho:\C\times
\C^5\to \C^5$ denote the
principal right action on the $\C$-principal bundle $\pi:\C^5\to X$
and let $Y=\pi^{-1}(H_0)$. We define a $\C$-action on $Y$ by
\[
\mu_f(t):y\mapsto \rho(t\cdot f(\pi(y)),y)
\]
and claim that for $t\ne 0$ the automorphism $\mu_f(t)$ of the
hypersurface $Y$ can not be extended to an automorphism of $\C^5$,
neither algebraically nor holomorphically.

Fix some $t\ne 0$ and assume that $\phi\in \Auth(\C^5)$ is an extension
of $\mu_f(t)$. There is no loss in generality in assuming $t=1$.
Then $\phi$ stabilizes $\C^5\setminus Y$ and therefore
induces a holomorphic map $\pi\circ\phi:\C^5\setminus Y\to X\setminus H_0$.
We recall that $X\setminus H_0$ embeds into $\Qbar$ in such a way
that the complement is $D\cup H$, $D\cup H$ is a divisor with simple
normal crossings and $D+H+K$ is ample on $\Qbar$.
By a theorem of Griffiths and King (\cite{GK}, Prop.~8.8) this implies that
$\pi\circ\phi$ is algebraic. Furthermore $X\setminus H_0$ is of log
general type as defined by Iitaka \cite {I}.
We claim that $\phi$ must map $\pi$-fibers into $\pi$-fibers.
Indeed, otherwise there would exist an irreducible algebraic
subvariety $R$ of codimension two in $\A^5$ such that there is a
dominant morphism $F:R\times\C\to X\setminus H_0$ given by
\[
F(r,t)=\pi\circ\phi(\rho(t,r))
\]
Since $X\setminus H_0$ is of log general type, such a map can not exist.
It follows that $\pi$ is equivariant for the automorphism $\phi$.
However, being of log general type $X\setminus H_0$ admits only
finitely many automorphisms \cite{I},\cite{Sa}.
Therefore there is a number $m$ such that
$\phi^m$ induces the trivial action on the base, i.e.~$\phi^m$ acts
only along the $\pi$-fibers.

Let $(U_i)_{i\in K}$  denote a covering of $X$ by open affine
subsets.
The $\C$-principal bundle may be described in terms
of transition functions $\zeta_{ij}\in\C[U_i\cap U_j]$.
With respect to the corresponding local trivialization $\psi=\phi^m$
is given by $(p,x)\mapsto (p,\alpha_i(p)x +\beta_i(p))$
for $(p,x)\in U_i\times\A^1$ with $\alpha_i\in\O^*(U_i)$ and
$\beta_i\in\O(U_i)$.
An easy calculation shows that $\alpha_i=\alpha_j$ and
$\beta_i=\beta_j + (\alpha_i-1)\zeta_{ij}$ on $U_i\cap U_j$.
Thus $\alpha=\alpha_i$ is a global holomorphic function defined on the
whole $X$. Since $\codim(E)=2$, the function $\alpha$ can actually
be defined on the affine variety
$Q^*=Q\setminus S$. Recall that the transition functions $\zeta_{ij}$ are
algebraic functions on $U_{ij}$
 and therefore extend to rational functions on the whole $Q$.
It follows that $(\alpha-1)\zeta_{ij}$ can be extended to a meromorphic
function on $Q$ for all $i,j$. Using $\beta_i=\beta_j
+(\alpha-1)\zeta_{ij}$, this implies that the functions $\beta_i$ can
be extended to meromorphic functions on $Q^*=\cup_i U_i$.
Moreover, if $h$ is a regular function on $Q^*$ such that
none of the functions $h\zeta_{ij}$
has poles on $Q^*$, then the functions
$h\beta_i$ haven't any poles in
$Q^*$ either.

Now observe that $\alpha|_{H_0\cap Q^*}\equiv 1$ and
$\beta_i|_{H_0\cap Q^*}\equiv mf$ for all $i$, because $\psi=\phi^m$ coincides
with
$\mu_f(m)$ on $\pi^{-1}(H_0\cap Q^*)$.
It follows that, for $h$ chosen as above,
 $hf$ defines a holomorphic
function on $(H_0\cap Q^*) \cup E$.
However, the condition which $h$ has to fulfill
(namely, that all the $h\zeta_{ij}$ are regular on $Q^*$) does not
depend on the choice of $f$.
Hence it is possible to choose $f$ and $h$ in such a way that
$hf$ has poles in $E$ and in this case the assumption of
the existence of $\phi$ leads to a contradiction.

Thus it is possible to choose $f$ such that the resulting group
action $\mu_f(t)$ on the hypersurface $Y$ in $\C^5$ can not be extended
to an algebraic or holomorphic automorphism of $\C^5$.
\end{proof}
\begin{remark} If one is interested only in the algebraic
non-extendibility, then the last steps of the proof can be simplified
substantially:

If $\alpha$ is algebraic, it is necessarily constant.
This implies $\alpha\equiv 1$, since $\alpha|_{H_0\cap Q^*}\equiv 1$.
Hence
$\beta=\beta_i$ is a global regular function and the non-extendibility
of $\mu_f(t)$ follows directly from the fact that $f$ does not extend
to a regular function on $X$.
\end{remark}
\begin{remark}
Every $\C$-principal bundle over a differentiable manifold
is differentiably trivial.
Using this fact, it is easy to see that the $\C$-action on $H$
does extend to a differentiable action on $\C^5$.
\end{remark}

\section{Extending actions of whole groups}
As  mentioned above, every affine variety $Z$ may be embedded
into some affine space $\A^n$ in such a way that every automorphism
of $Z$ extends to an automorphism of $\A^n$.
If $\Aut(\A^n,Z)$ denotes the group of all automorphisms of $\A^n$
stabilizing $Z$ (as a set, not pointwise), this is equivalent to the
assertion that the natural group homomorphism
$\Aut(\A^n,Z) \to \Aut(Z)$
is surjective.
Thus there exists a short exact sequence
\[
1 \to L \to \Aut(\A^n,Z) \to \Aut(Z) \to 1
\]
one can ask if it splits.

Here we will show that there is an affine variety,
namely $Z=\A_1^*\times\A_1^*$, such that this sequence never splits.
More precisely we will prove the following theorem.
\begin{theorem}
Let $k$ be an algebraically closed field,
$K$ a field, $n\in\N$
and let $G$ denote the group given by the
semidirect product $SL_2(\Z)\sdr_\rho (k^*\times k^*)$
with
\[
\rho\begin{pmatrix} a & b \\ c & d \\ \end{pmatrix}(z,w)=
(z^aw^b, z^cw^d)
\]

Then there does not exist any injective group homomorphism from
$G$  either into the group $\Diff(\RR^n)$ of diffeomorphisms of $\RR^n$
or into the group $\Aut_K(\A^n)$ of $K$-automorphisms of the affine
space $\A^n$.
\end{theorem}
\begin{corollary}
Let
either
$Z$ denote the complex manifold $\C^*\times\C^*$ and $H$ be the
group of holomorphic automorphisms of $Z$ or let $k$ be an
algebraically closed field, $Z=\A^2\setminus\{(z_1,z_2):z_1z_2=0\}$
and $H$ denote the group of $k$-automorphisms of $Z$.

Then there does not exist any non-constant
$G$-equivariant holomorphic resp.~regular map from $Z$ to
$\C^n$ resp. $\A^n$.
\end{corollary}
\begin{proof}
By the theorem there is no injective group
homomorphism from $G$ into $\Auth(\C^n)$ resp.~$\Aut_k(\A^n)$.
Since $G\subset H$, it follows that
an equivariant map $F$ must have non-trivial fibers such that
a non-trivial normal subgroup of $H$ acts only in fiber direction.
Now $Z$ is homogeneous under the $H$-action and
it is easy to check that the only $H$-equivariant fibrations of $Z$
are the maps $\tau_N:Z\to Z$ given by $(z_1,z_2):
\mapsto (z_1^N,z_2^N)$.
Since the base of these fibrations is again isomorphic to $Z$,
it follows from the theorem that there is no non-constant
equivariant holomorphic map
resp.~morphism from $Z$ to $\C^n$ resp.~$\A^n$.
\end{proof}

As a first step towards  theorem~\ref{harmful-theo} we need
the following well-known result on the
existence of fixed points.

\begin{proposition}\label{fix-pt}
Let $p$ be a prime and let
$\Gamma$ be a finite abelian $p$-group acting differentiably on $\R^n$ or
algebraically on an affine space $\A^n$ defined over a field $k$
with $char(k)\ne p$.

Then $\Gamma$ has a fixed point.
\end{proposition}
\begin{proof}
In the first case this follows from Smith theory \cite{H},
in the second case a proof may be sketched as follows:
Every $\gamma\in\Gamma$ induces an automorphism of $\A^n$ given
by a polynomial map.
Let $R\subset k$ be the ring generated by all the
coefficients of these polynomials for all $\gamma\in\Gamma$.
Let $\p$ be a prime ideal in $R$
and consider reduction modulo $\p$. Then $R/\p$ is a finite field.
If $char(R/{\p})\ne p$, then $p$ does not divide the number of
points in $\A^n(R/\p)$. This implies that there must be a fixed point
modulo $\p$. Finally, the existence of fixed points for almost all
prime ideals in $R$ implies that there is a fixed point in $\A^n(k)$.
\end{proof}

\begin{lemma}\label{eff-tang}
Let $p$ be a prime and
let $\Gamma$ be a finite $p$-group.
Let $V$ be a connected differentiable manifold or an
irreducible variety defined
over a field $k$ with $p\ne char(k)$,
assume that $\Gamma$ acts effectively on $X$
and let $v$ be a fixed point.

Then $\Gamma$ acts effectively on the (Zariski-)
tangent space $T_vV$.
\end{lemma}
\begin{proof}
In the differentiable case finiteness of $\Gamma$ permits us to
construct an invariant Riemannian metric. Then $\Gamma$ acts by
isometries. Hence if $\gamma\in\Gamma$ fixes a point $v\in V$ and acts
trivially on the tangent space, it preserves all the geodesics
emanating from $v$ and therefore induces the identity map on $V$.

In the algebraic case we note that, if $\gamma\in\Gamma$ fixes $v$
and acts trivially on the Zariski tangent space
$T_vV=(m_v/m_v^2)^*$, then $\gamma$ acts trivially on $m_v^k/m_v^{k+1}$
for all $k\in\NN$. Let $f\in m_v$ with $\gamma^*(f)\ne f$.
Then $\gamma^*(f)-f\in m_v^k\setminus m_v^{k+1}$ for some $k\in\NN$.
But this implies that
$(\gamma^*)^n(f)-f \equiv n \left( \gamma^*f-f \right)$
modulo $m_v^{k+1}$. This contradicts our assumption that $\gamma$ is
of finite order coprime to the characteristic of $k$.
Hence the assertion.
\end{proof}

Now we are in a position to prove the theorem.

\begin{proof}
If there would be such an injective group homomorphism,
we would obtain an effective $G$-action on $\R^n$ resp.~$\A^n$.

Fix a prime $p$, $p\ne char(k),char(K)$,
let $A=\C^*\times\C^*$ resp.~$A=k^*\times k^*$ and let
$\Gamma_r=A[p^r]$ denote the subgroup of torsion elements of order
$p^r$ in $A$.
Let $N_G(\Gamma_r)$ resp.~$Z_G(\Gamma_r)$ denote the normalizer
resp.~centralizer of $\Gamma_r$ in $G$.
Note that $\Gamma_r\simeq (\Z/p^r\Z)\times(\Z/p^rZ)$,
$N_G(\Gamma_r)/Z_G(\Gamma_r)\simeq SL_2(\Z/p^r\Z)$
and that $N_G(\Gamma_r)/Z_G(\Gamma_r)$ acts effectively on $\Gamma_r$
by group automorphisms.

Let $V_r$ denote the fixed point set of the
induced $\Gamma_r$-action on
$\R^n$ resp.~$\A^n$.
We claim that there exists some uniform (i.e.~independent of $r$)
upper bound $C<\infty$ for the number of (Zariski-)connected components of
$V_r$.
In the differentiable case it follows from Smith theory (see e.g.~\cite{H})
that $V_r$ is connected.
In the algebraic case the sequence $V_r$ is a descending sequence
of subvarieties of $\A^n$ and  hence there is some $R$ such that $V_r=V_R$
for all $r\ge R$. This yields the desired bound $C$ for the number
of Zariski connected components of $V_r$.

On the other hand, proposition \ref{fix-pt}
 implies that none of the $V_r$ can be
empty.
Lemma \ref{eff-tang} implies that $\Gamma_r$ acts effectively on $T_x$ for
every $x\in V_r$.
The type of the representation of $\Gamma_r$ on $T_x$ is determined
by the trace function $tr_x$ defined by $\gamma\mapsto
Tr(\rho_x(\gamma))$ where $\rho_x$ denotes the action of $\Gamma_r$
on $T_x$. This trace function must depend continuously resp.~regularly
on $x$. However, there are only finitely many possible values.
Hence, the trace function is locally constant with respect to $x$,
i.e., on each (Zariski-) connected component of $V_r$ all the
representations $\rho_x$ of $\Gamma_r$ are isomorphic.

For every number $r\in\N$ fix a point $x_r\in V_r$ and let
$N^0_r$ denote the subgroup of those elements of $N_G(\Gamma_r)$ which
stabilize the connected component of $V_r$ containing $x_r$.
For every $\phi\in N^0_r$, $\gamma\in\Gamma_r$, the representation
of $\gamma^\phi=\phi\gamma\phi^{-1}$ and $\gamma$ on $T_{x_r}$ are
isomorphic, because $\phi(x_r)$ and $x_r$ lie in the same connected
component of $V_r$. It follows that the representations differ only
by a permutation of the irreducible $\Gamma_r$-submodules of
$T_{x_r}$.
It follows that $\# N^0_r/Z_G(\Gamma_r)\le n^!$.
Thus $\#N(\Gamma_r)/Z(\Gamma_r) \le C (n^!)^2$.

But $\#N(\Gamma_r)/Z(\Gamma_r)$ goes to infinity for increasing $r$,
because $N(\Gamma_r)/Z(\Gamma_r)\simeq SL_2(\Z/p^r\Z)$.
Thus we deduced a contradiction from the assumption that there exists
such an injective group homomorphism.
\end{proof}

\section{The Analytic Case}

Here and in the rest of our paper holomorphic embedding means a proper
holomorphic embedding, i.e., a proper holomorphic map which is
injective and immersive.

We will use
the ideas developed in \cite{RR1} together with the techniques for
interpolation by
automorphisms of $\C^n$ developed in \cite{AL}, \cite{BF}, \cite{FGR},
\cite{F1}, \cite{FR} to construct
our more complicated embeddings for a complex subspace $X$ of $\C^n$ of any
dimension.

If $X$ is a Stein manifold of dimension $n$, then
$X$ can be embedded into $\C^N$ with $N=\max\{3,n+[\frac{n}{2}]+1\}$.
This bound is optimal, except possibly for $n=1$.
(see \cite{Sch}).

We want to point out that we give no results about existence of embeddings
of a given Stein space into
$\C^n$.%

Instead we consider the question which further properties of such
an embedding  can be prescribed,
if there already exists an embedding in the given dimension.
Further results
in this direction can be found in \cite{F1}, \cite{BF}, \cite{BFo}, \cite{RR2}.

\subsection{Results}

Let $T$ be a topological space.
By a continuous family of holomorphic selfmaps of a complex space $X$
parameterized by $T$ we mean a continuous map $\alpha : T\times X \to X$ such
that for every fixed $t\in T$ the map $\alpha (t,\cdot): X\to X$ is
holomorphic.
If $T$ is a real Lie group such that $\alpha$ describes a group
action, then the map $\alpha$ is automatically real
analytic (see for instance \cite{A}).
We will call a holomorphic map $F:\C^n
\to \C^n$  nondegenerate if the Jacobian $JF$ of $F$ does not
vanish identically.  Our main result is the following theorem,
which will be proved in the next paragraph.

\addtocounter{theorem}{1}
\begin{theorem}\label{main-ana}
Let $X$ be a proper complex analytic subvariety of $\C^n$ which consists
of infinitely many points and $\alpha :T\times X
\to X$ a continuous family of holomorphic selfmaps of $X$ parametrized
by a locally compact topological space $T$ with countable topology. Then there
exists an embedding $\varphi: X\ha \C^n$ with the following property:
\smallskip
\noindent
If $F:\C^n\to\C^n$ is a nondegenerate holomorphic map such that
\begin{enumerate}
\item  $F^{-1}(\C^n\setminus \varphi (X))= \C^n\setminus \varphi(X)$ and
\item  $\varphi^{-1}\circ F\circ\varphi=\alpha(t,\cdot)$ for some
$t\in T$,
\end{enumerate}
then $F= id$.
\end{theorem}

Theorem 3 stated in the introduction then arises as a corollary.

\begin{corollary}
Let $X$ be a complex space of dimension $m$ which can be embedded into $\C^n$
and such that the holomorphic automorphism group $\Auth (X)$ is a Lie group.
Also assume that $m<n$ and $X$ is not a finite set of points.
Then there exists an embedding $\varphi: X\ha \C^n$ of $X$ into
$\C^n$ such that the only holomorphic automorphism of $\C^n$ leaving
$\varphi (X)$ invariant is the identity.
If furthermore $n-\dim (X)\ge 2$, then
$\phi$ can be chosen so that
the holomorphic automorphism group $\Auth (\C^n\setminus \varphi (X))$ of the
complement is trivial.
\end{corollary}
\begin{proof}
By definition, a Lie group is locally compact and of countable topology.
Hence we can apply the theorem to $X$
with $T=\Auth (X)$.
 Let $\varphi: X\ha \C^n$ be the embedding given by theorem~\ref{main-ana}.
Any holomorphic automorphism $\tilde \alpha $ of $\C^n$ leaving
$\varphi (X)$ invariant
is clearly a nondegenerate holomorphic map from $\C^n$ to
$\C^n$ with $\tilde \alpha^{-1} (\C^n \setminus \varphi (X))=
\C^n\setminus\varphi (X)$ and the restriction of $\tilde\alpha$ to
$\varphi (X)$
is an holomorphic automorphism of $X$, hence an element of our family. This
imply that $\tilde\alpha$ is the identity map.

For the last statement of the corollary observe that, since
$\varphi (X)$ has at least codimension 2 in $\C^n$, any holomorphic
automorphism $\alpha \in \Auth (\C^n\setminus \varphi(X))$ of the complement
extends to an holomorphic automorphism $\tilde \alpha$ of $\C^n$ and
which leaves
$\varphi (X)$ invariant.
\end{proof}
\begin{corollary}
Let $X$ be a Stein manifold and let $G$ be a Lie group acting
effectively and by biholomorphic transformations on $X$.

Then there exists a natural number $N$ and an embedding
$X\hookrightarrow \C^N$ such that for no $g\in G\setminus\{e\}$ the
induced automorphism of $X$ can be extended to an automorphism
of $\C^N$.
\end{corollary}

If $k=1$, then $\Auth(\C^k)$ is a Lie group.
Hence we obtain the following corollary.

\begin{corollary}
For every $n\ge 2$ there exists an embedding $\varphi: \C \ha \C^n$
 such that the only holomorphic automorphism of $\C^n$
leaving $\varphi (\C)$ invariant is the identity. If $n\ge 3$,
 then $\phi$ can be chosen such that
 $\Auth (\C^n\setminus \varphi (\C))$
 is trivial.
\end{corollary}

\begin{question}
Does there exist an embedding $\varphi:\C^k\to\C^n$ with
$\Auth(\C^n\setminus \varphi(\C^k))=\{id\}$ for every $0<k<n$ ?
\end{question}

\begin{proposition}
Let $X\subset\C^n$ be a complex analytic hypersurface.

Then there exists an embedding $\varphi:X\to\C^n$ such that
$\Auth(\C^n\setminus X)$ is countable.
\end{proposition}
\begin{proof}
In \cite{BF} Buzzard and Fornaess  proved (see also \cite{F1} Theorem 5.1  for
a more
general statement) that there is an embedding $\varphi : X\ha \C^n$ such that
the complement $\C^n\setminus \varphi (X)$ is Kobayashi-hyperbolic.
Hence the result follows from the lemma below.
\end{proof}
\begin{lemma}
Let $V$ be an affine algebraic manifold (e.g.~$V=\C^n$) and
$Y\subset V$ be a complex analytic subvariety such that the complement
$V\setminus Y$ is Kobayashi-hyperbolic.
Then  $\Auth (V\setminus Y))$ is countable and discrete in the
compact-open topology.
\end{lemma}

\begin{proof}
The group of holomorphic automorphisms of a Kobayashi-hyperbolic manifold
is a closed subgroup of the isometry group with respect to the Kobayashi
distance and is isomorphic to  a real Lie group.
We must show that the identity component of
$\Auth (V\setminus Y)$ is zero-dimensional.
If $\Auth (V\setminus Y)$ were positive dimensional,
then any nonzero vector $v$ in the Lie algebra would give rise to
a non-trivial action of the additive group $(\R,+)$ on $V\setminus Y$
given by the map $(t,z)\mapsto exp (tv)\cdot z$. Since $V\setminus Y$ is
Stein and has no non-constant bounded plurisubharmonic functions  we can
apply Cor. 2.2.
in \cite{F} to conclude that our action of $(\R,+)$ extends to an
$(\C,+)$ action on $V\setminus Y$. This contradicts the fact that every
holomorphic map from $\C$ into the
Kobayashi-hyperbolic manifold $V\setminus Y$ is
constant.
\end{proof}

\begin{question}
Given a complex analytic hypersurface $X\subset\C^n$, does there
exists an embedding $j:X\hookrightarrow \C^n$ such that
the automorphism group of $\C^n\setminus j(X)$ is trivial ?
\end{question}

More generally, one may pose the following question.

\begin{question}
Given a complex analytic hypersurface $X\subset \C^n$ and a connected
complex Lie group $G$, under which conditions it is true that there
exists an embedding $j:X\hookrightarrow \C^n$ such that
$\Auth(\C^n\setminus X)=G$?
\end{question}
One may also ask the same question for
$\{\phi\in\Auth(\C^n):\phi(j(X))=j(X)\}$ instead
of $\Auth(\C^n\setminus(j(X)))$.
\subsection{Proofs}

As already mentioned above, in the proof of Theorem~\ref{main-ana}
we will use the techniques
developed in \cite{F1}, \cite{RR1} (see also \cite{FGR}, \cite{BF}).
For the convenience of the reader
we state here all results from these papers which will be used in the sequel.

\begin{proposition}\label{prop:f1}
(see \cite{F1} Prop. 1.1.)  Let $n>1$. Assume that

\begin{itemize}
\renewcommand{\labelitemi}{\addtocounter{enumi}{1}\alph{enumi})}
\setcounter{enumi}{0}
\item $K\subset \C^n$ is a compact polynomially convex set,

\item $\{a_j\}_{j=1}^s\subset K$ is a finite subset of $K$,

\item $p$ and $q$ are arbitrary points in $\C^n\setminus K$ (not necessarily
distinct)

\item $N$ is a nonnegative integer, and

\item $P:\C^n\to\C^n$ is a holomorphic polynomial map of degree at
   most $m\ge 1$ with $P(0)=0$ and $JP (0)\ne 0$.
\end{itemize}

Then for each $\epsilon >0$ there exists an automorphism $F\in \Auth (\C^n)$
satisfying

\begin{itemize}
\renewcommand{\labelitemi}{\addtocounter{enumi}{1}\roman{enumi})}
\setcounter{enumi}{0}
\item $F(p)=q$ and $F(z)=q +P(z-p)+ O(\vert z-p\vert^{m+1})$ as $z\to p$,

\item $F(z)=z+O(\vert z-a_j\vert^n )$ as $z\to a_j$ for each $j=1,2,\ldots,s$
and

\item $\vert F(z) -z \vert + \vert F^{-1} (z) -z\vert < \epsilon \quad \forall
\ z \in K$
\end{itemize}

\end{proposition}

\begin{proposition} \label{prop:f2}
(see \cite{F1} Prop. 4.1. and 4.2.)
Let $K_0\subset K_1\subset K_2\subset \dots\subset \cup_{j=0}^\infty K_j =\C^n$
be compact sets such that $K_{j-1}\subset {\rm Int}K_j$ for each $j\in\N$.
Suppose $\epsilon_j$ $(j=1,2,3,\ldots)$ are real numbers such that
$$0<\epsilon_j<\dist (K_{j-1},\C^n\setminus K_j)\ \ (j\in\N), \quad
\sum_{j=1}^\infty \epsilon_j <\infty . \hfill\eqno{(1)}$$
Suppose that for each $j=1,2,3,\ldots$ $\Psi_j$ is a holomorphic automorphism
of $\C^n$ satisfying
$$\vert \Psi_j (z) -z \vert <\epsilon_j,\quad z\in K_j. \hfill \eqno{(2)}$$
Set $\Phi_m=\Psi_m\circ\Psi_{m-1}\circ\dots\circ\Psi_1$. Then there is an open
set $\Omega\subset \C^n$ such that $\lim_{m\to\infty} \Phi_m =$ $\Phi$
exists on $\Omega$ (uniformly on compacts), and $\Phi$ is a biholomorphic
map of $\Omega$ onto $\C^n$. In fact, $\Omega =\cup_{j=1}^\infty \Phi_m^{-1}
 (K_m)$. Also $\Omega$ can be characterized as the set of points $z\in\C^n$
such that the sequence $\{ \Phi_m (z): m\in\N \}$ is bounded.
\end{proposition}

As usual $B_{r}$ denotes the (open) ball of radius $r>0$ in $\C^n$
and $\partial B_r$ denotes its boundary.

\begin{lemma}  (see \cite{RR1} Lemma 4.3.)
Let $0<a_1<a_2$, $0<r_1<r_2$, $c>0$ be real numbers and $A\subset
\partial B_{r_1}$ a dense subset.
Then there exists a finite subset
$\bigcup_{i=1}^m z_i$ of $A$ such that
$$
F(B_{a_1})\subset B_{r_1}
$$
 for all holomorphic maps
$F: B_{a_2}\to B_{r_2}$
with
$F(0)\in B_{\frac{{r_1}}{2}}$,
$|JF (0)|\ge c$
and $z_i\notin F( B_{a_2})\quad \forall \ i=1,2,\ldots, m$.

\end{lemma}
\begin{remark}
Lemma 1 is proved in \cite{RR1} with $A=\partial \B_{a_1}$. The
proof starts with an arbitrary countable dense subset $x_1,x_2,\ldots$
of $\partial \B_{a_1}$ and the desired finite set is constructed as a
subset of $x_1,x_2,\ldots$. So the only minor modification to be made in
the proof is to start with $x_1,x_2,\ldots$ being a subset of $A$ which
is dense in $\partial \B_{a_1}$.
\end{remark}

The proof of theorem~\ref{main-ana} consists of two steps.
The first step  is to
construct an  embedding in such a way
that conditions 1.\ and 2.\ of the theorem force
a  nondegenerate holomorphic map $F:\C^n \to \C^n$ to be affine, i.e., an
affine automorphism of $\C^n$. The second point is to ensure that the only
affine automorphism of $\C^n$ leaving the embedded variety $X$ invariant is
the identity. To be accurate in the second point we need one more technical
result to be explained now:

By a submanifold $Z\subset \C^n$ we shall mean an injectively immersed
holomorphic submanifold (not necessarily closed).

\begin{definition}
Let $k\ge 2$ be a natural number. We say that a submanifold $X$
of $\C^n$ osculates of order $k$ at some point $x\in X$ if $X$
has contact order $k$ with the tangent space $T_xX\subset \C^n$ at
$n$.
\end{definition}
In local coordinates this condition can be expressed in the following
way:

Let $\varphi :U(\subset \C^m) \to X$ be a holomorphic coordinate chart
for the
$m$-dimensional manifold $X$
around $x$ (i.e. $\varphi (0)=x$).
Then $X$ osculates of order $k$ at $x$ if and only if
$\frac{\partial }{ \partial w^\alpha }
\vert_{w=0} \varphi\in T_x X$ for every multiindex
$\alpha= (\alpha_1,\alpha_2,\ldots,\alpha_m)$ with
$2\le\vert \alpha\vert\le k$.

The property to osculate of order $k$ is evidently preserved by
affine coordinate changes on $\C^n$, i.e., if $\psi :\C^n\to \C^n$ is an
affine automorphism of $\C^n$, then a submanifold $X\subset \C^n$ osculates
of order $k$ at $x\in X$ iff the submanifold $\psi (X)$ osculates of
$k$ at $\psi (x)\in \psi (X)$.

\begin{lemma}\label{klem:2}
Let $M\subset \C^n$ be an $m$-dimensional $(m<n)$ submanifold.
Suppose we are given

\begin{itemize}
\setcounter{enumi}{0}
\renewcommand{\labelitemi}{\addtocounter{enumi}{1}\alph{enumi})}
\item  a compact subset $K_M$ of $M$,

\item  a compact subset $K$ of $\C^n$,

\item  finitely many points $a_1,a_2,\ldots,a_r$ in $K_M\subset M$,

\item  finitely many points $b_1,b_2,\ldots,b_q$ in $M\setminus K_M$,

\item a natural number $k\ge 2$, if $m=1$ and $n=2$ then $k\ge 3$, and

\item a real number $\epsilon >0$.
\end{itemize}

Then there exists an holomorphic automorphism $\psi\in \Auth (\C^n)$ with the
following properties:
\begin{enumerate}

\item $\psi (a_i)=a_i\quad \forall\ i=1,2,\ldots,r$

\item $\psi (z)=z+O (\vert z-b_i\vert^{k+1})$ as $z \to b_i$
      $\forall\ i=1,2,\ldots,q$

\item $\vert \psi(z)-z \vert + \vert \psi^{-1} (z) -z\vert <\epsilon\ \forall\
z\in K$

\item there is no point $m\in K_M$
such that the submanifold $\psi (M)$ of $\C^n$ osculates
of order $k$ at $\psi (m)$.

\end{enumerate}

\end{lemma}

Before we can prove lemma~\ref{klem:2} we need a sublemma.
The notations used in
the sublemma are the same as those in lemma~\ref{klem:2}.

\begin{sublemma}
For each point $m\in K_M$ there is an open neighborhood $U_m$ of
$m$ in $M$ and a family $\psi_t$ of
automorphisms of $\C^n$ parametrized by $\C^{\mk}$
with $\mk=\mkn$ such that
\begin{enumerate}
\item
$\psi_0=id$,
\item
Every $\psi_t$ fulfills conditions $(1)$ and $(2)$ of
lemma~\ref{klem:2}
\item
There exists an open neighbourhood $T$ of $0$ in $\C^{\mk}$
such that
\[
\Sigma=\{t\in T:\exists m'\in U_m:
\text{ $\psi_t (M)$ of $\C^n$  osculates of order $k$ at
$\psi_t (m^\prime )$}\}
\]
is a set of Lebesgue measure zero.
\end{enumerate}

\end{sublemma}

\begin{proof}
If $M$ does not osculate of order $k$ at $m$,
then $\psi_t=\id$ for all $t$
does the job.
Hence we may assume that $M$ osculates of order $k$ at $m$.

Without loss of generality we may assume that $m=0$,
$T_mM=\{(z_1,\ldots,z_m,0,\ldots,0\}$.
Let $\pi:\C^n\to\C^m$ denote the map given by projection onto the
first $m$ coordinates. After a linear change of coordinates we may assume that
$\pi(b_i)\ne 0\in\C^m$ for all $1\le i\le q$.

To examine whether $\psi(M)$ osculates of order $k$ at some point
$\psi(m^\prime)$ with $m'\in M$ for a given $\psi\in\Auth$
we consider the map
$F^\psi: M \to \C^\mk$ whose coordinate
functions $F^\psi_{\alpha, u}$ are enumerated by
pairs $(\alpha,u)$ where $\alpha$ is a multiindex
$\alpha=(\alpha_1,
\alpha_2,\ldots,\alpha_m)$  with $2\le \vert \alpha\vert\le k$ and
$u\in\NN$ satisfies
 $m+1\le u\le n$ ,
these are given by
$$
F^\psi_{\alpha,u} (w) = \det \left(
\begin{array}{cccc}
\dd{w_1} (\psi)_1 (w) &\ldots &
\dd{w_1} (\psi)_m (w) &
\dd{w_1} (\psi)_u (w) \\
\vdots& &\vdots&\vdots \\
\dd{w_m} (\psi)_1 (w) &
\ldots &
\dd{w_m} (\psi)_m (w) &
\dd{w_m} (\psi)_u (w) \\
\dd{w^\alpha} (\psi)_1 (w) &
\ldots &
\dd{w^\alpha} (\psi)_m (w) &
\dd{w^\alpha} (\psi)_u (w)
\end{array}
\right).
$$
Here $(\psi_t)_i$ denotes the $i$-th coordinate function of the
map $\psi:\C^n\to\C^n$ and $(w_i)_{1\le i\le m}$ is some
fixed system of local coordinates on $M$ near $m$.
Then $\psi(M)$ osculates of order $k$ at $w$ if and only if
$F^\psi(w)=0$.

By restricting our attention to a small enough neighbourhood of
$m\in M$ we may choose $w_i=z_i$ ($1\le i\le m$).

Now we come to the explicit construction of our family $(\psi_t)_t$ of
automorphisms of $\C^n$.

For each pair $(\alpha, u)$ with $\alpha$ and
$u$ as above we choose a holomorphic function $h_{\alpha,u}$
on $\C^m$ such that
\begin{enumerate}
\item
$h_{\alpha,u}-z^\alpha$ vanishes of order at least $k+1$ in $0$.
\item
$h_{\alpha,u}$ vanishes of order at least $k+1$ for all
$\pi(b_i)$ ($1\le i\le q$),
\item
$h_{\alpha,u}$ vanishes at all $\pi(a_j)$ ($1\le j\le r$).
\end{enumerate}

Next we  define
a map $\psi:\C^\mk\times\C^n\to\C^n$ by
\[
\psi(t,z)= z + \sum_{(\alpha,u)} t_{(\alpha,u)}
h_{\alpha,u}(z_1,\ldots,z_m)e_{u}
\]
Here $e_u$ denotes the $u$-th unit vector and the coordinates of
$\C^\mk$ are indexed by pairs $(\alpha, u)$ as above.
For every $t\in\C^{\mk}$ the map $\psi_t=\psi(t,\cdot)$
is an automorphism of $\C^n$ fulfilling (because of 2. and 3.)
conditions
$(1)$ and $(2)$ of lemma~\ref{klem:2}.
Furthermore $\psi_0=id$.

Easy calculations (using 1.) show that
\[
\dd{t_{(\alpha,u)}}F^{\psi}_{\alpha,u}|_{w=0}
= \alpha !
\]
and
\[
\dd{t_{(\alpha^\prime,u^\prime)}}F^{\psi}_{\alpha,u}|_{w=0} = 0
\]
whenever $u\ne u^\prime$ or whenever $u=u^\prime$, $|\alpha^\prime|\le|\alpha|$
and $\alpha^\prime\ne \alpha$.

This implies that the map $\Phi:\C^{\mk}\times M\to\C^{\mk}$
defined by $\Phi(t,z)= F^{\phi_t} (z)$ has maximal rank near $0$.
Thus there exists an open neighborhood $\Omega_m$
of the form $\Omega_m=T\times U_m$ of $(0,m)$ in $\C^\mk\times M$ such that
$\Phi |_{\Omega_m}$ is transversal to $0\in \C^\mk$.
This implies that for
almost all $t\in T$ the map $F^{\psi_t} : U_m \to \C^\mk$ is transversal
to $0$. Since $m<\mk$, (Here we need $k>2$ if $n=2$ and $m=1$. In all other
cases $k=2$ is already sufficient.) this means that for almost all $t\in T$
the image $F^{\psi_t} (U_m)$ does not meet $0$,
i.e.~$\psi_t(M)$ does not osculate of order $k$ for any $m'\in U_m$.
\end{proof}

\begin{proof}[Proof of Lemma~\ref{klem:2}.]
Choose finitely many open subsets $U_i$ of $M$ together with families
$\psi^i: T_i\times\C^n\to \C^n$ of automorphisms  $i=1,2,\ldots, l$
as in the sublemma and choose
 compact subsets $K_i\subset U_i$ of the $U_i$ which cover $K_M$.
Since $\psi^1_0$
is the identity, for $t$ sufficiently small the automorphism $\psi^1_t$ moves
no point of $K_M$ more then $\frac{\epsilon}{l}$. So we find a $t_1\in T_1$
such that $\vert \psi^1_{t_1} (z) - z \vert < \frac{\epsilon}{ l }$ $\forall \
z\in K$
and the submanifold $\psi^1_{t_1} (M)$ does not osculate of order
$k$ at any point of $\psi^1_{t_1} (K_1)$.

Observe that the property of
osculating of order $k$ at some point is preserved under small perturbations,
i.e., for each compact subset $L$ of a submanifold $M$ of $\C^n$ which
does not osculate of order $k$ at any point of
$L$ there exists some $\epsilon$ such
that for each automorphism $\Psi$ of $\C^n$ the property
$\vert \Psi (z) -z \vert <\epsilon \quad \forall \ z\in L$  implies that
$\Psi (M)$ remains non-osculating of order $k$ at any point of $\Psi (L)$
(For holomorphic maps small perturbations in values imply small perturbations
in derivatives). Hence we find a sufficiently small $t_2 \in T_2$ such that
first $\vert \psi^2_{t_2} (z) - z \vert < \frac{\epsilon}{ l }$
$\forall \ z\in \psi^1(K)$, second the submanifold
$\psi^2_{t_2} \circ\psi^1_{t_1} (M)$
does not osculate of order $k$ at any point of
$\psi^2_{t_2} \circ\psi^1_{t_1} (K_2)$ and third
$\psi^2_{t_2} \circ\psi^1_{t_1} (M)$ remains non osculating of order $k$ at
any point of $\psi^2_{t_2} \circ\psi^1_{t_1} (K_1)$. Proceeding by induction
we find an automorphism $\psi:=\psi^l_{t_l} \circ \psi^{l-1}_{t_{l-1}}\circ
\dots\circ \psi^1_{t_1}$ moving no point of $K$ more than $\epsilon$ and
such that $\psi (M)$ does not osculate of order $k$ at any point of
$\psi (\cup_{i=1}^l K_i) \supset \psi (K_M)$. Since all automorphisms
$\psi^i_t$ satisfy properties (1) and (2), $\psi$ satisfies them as well.
\end{proof}
\begin{proof}[Proof of Theorem 3.]
In the case that $X$ is a countable set of points the theorem is proved
in \cite{RR1} where these sets are called rigid. So we will deal only with the
case where $X$ is of positive dimension.

Let $\rho : X\to \R^{\ge 0}$ be a continuous exhaustion function, i.e.,
 $X_r:=\rho^{-1} ([0,r]) $ is a compact subset of $X$ for all $r\ge 0$.
Also let us denote by $\tilde X$ the union of all components of the smooth
part of $X$ which have maximal dimension, say $m$, $0<m<n$. We fix a natural
number $k\ge 2$. If $n=2$ and $m=1$ we require $k\ge 3$.
Let $b_1,b_2,\ldots,
b_{(n+2)}\in B_1\subset \C^n$ be points such that no affine automorphism of
$\C^n$ except the identity can permute these points.
Choose $x_1,x_2,\ldots,x_{n+2}\in \tilde X$ to be  mutually distinct points.
Applying Prop.~\ref{prop:f1}
several times (see also \cite{F1} Corollary 1.2.), we find
an automorphism $\Psi \in \Auth (\C^n)$
which maps $x_i$ to $b_i$ with the property that the submanifold $\Psi (\tilde
X)$
of $\C^n$ osculates of order $k$ at the points $\Psi (x_i)=b_i$ for each
$i=1,2,\ldots,n+2$. We denote the embedding $\Psi \circ i : X\ha  \C^n$
by $\f_0$.

Our aim is to construct a Fatou-Bieberbach domain
$\Omega \subset \C^n$ containing
$\f_0 (X)$ together with a biholomorphic map $\Phi : \Omega \to \C^n$ onto
$\C^n$ such that the embedding $\f:=\Phi\circ\f_0$ has the desired property.
The map $\Phi$ will be constructed as a limit of automorphisms $\Phi_m \in
\Auth (\C^n)$.

We start with an exhaustion $\cup_{i=1}^\infty T_i =T$ of the topological
space $T$ by compact subsets $T_i$ ($T_i\subset {\rm Int} T_{i+1}\subset T$).
Also we choose a sequence of open relatively compact
neighborhoods $U_i$ ($i=1,2,3,\ldots$) of the set $\cup_{i=1}^{n+2}\{
x_i\}$ in
$X$ with $\cap_{i=1}^\infty U_i = \cup_{i=1}^{n+2}\{ x_i\}$.

We will now  inductively define
real numbers
$\epsilon_m, R_m>0$,
natural numbers $k(m)$,
finite subsets $\cup_{j=1}^{k(m)}\{ a_j^m\}$
of $\partial B_{m+1}$, finite subsets $\{\cup_{j=1}^{k(m)}x_j^m\}$
of $X$, and automorphisms $\Phi_m$ of $\C^n$
for $m=0,1,\ldots$.
The beginning point is $\epsilon_0=\frac{3}{4}$, $R_o=1$, $k(0)=0$ and
$\Phi_0=\id$.
For $m\ge 1$  these data are recursively constructed in such a way
that the following
conditions are fulfilled:

\begin{itemize}

\item[($1_m$)]  $0<\epsilon_{m}<\frac{\epsilon_{m-1}}{3}$

\item[($2_m$)] If $F: \B_1 \to B_{m+2}\setminus \{\cup_{j=1}^{k(m)}
a_j^m\}$
is a holomorphic map with $\Vert F(0)\Vert \le \frac{m+1}{ 2}$ and
 $JF (0)\ge 1$ then
$F(B_{1-\frac{\epsilon_m}{ 2}})\subset B_{m+1}$.

\item[($3_m$)]  $\Phi_m \circ \f_0 (x_j^m) =a_j^m$ and $\rho (x_j^m)>
 \max_{x \in (\Phi_{m-1}\circ \f_0)^{-1} (\bar B_{m}) ,t\in T_m} \rho (\alpha
(t,x))$.

\item[($4_m$)] $\Vert \Phi_{m}\circ \f_0 (x) - \Phi_{m-1}\circ \f_0 (x)\Vert\le
\epsilon_{m}$
for all $x \in X_{ R_{m-1}}$

\item[($5_m$)] $\Vert \Phi_m \circ \Phi_{m-1}^{-1} (z) -z\Vert \le \epsilon_m$
for $z\in \bar B_m$

\item[($6_m$)] $\Phi_m\circ \f_0 (x_j^l)=a_j^l$ $\forall \ l<m\ j=1,2,\ldots,
k(l)$

\item[($7_m$)] $\Phi_m\circ\Phi_{m-1}^{-1} (z)=z +O(\vert z-b_i \vert^{k+1}) $
as $z\to
b_i$ for each $i=1,2,\ldots,n+2$

\item[($8_m$)] the submanifold $\Phi_m\circ \f_0 (\tilde X)$ of $\C^n$
does not osculate of order $k$ at any point $\Phi_m\circ\f_0 (x)$
with
$x \in (X_{ R_{m-1}}\cup \tilde X)\setminus U_m$

\item[($9_m$)] $\Vert \Phi_m\circ \f_0 (x)\Vert \ge m+1$
$\forall\ x\in X\setminus X_{R_m}$,

\item[($10_m$)] $R_m> R_{m-1}+1$

\end{itemize}

We will now confirm that such an recursive construction is
possible.
For  step 1 of the induction
we first choose $\epsilon_1 <\frac{\epsilon_0}{ 3}$
and, using Lemma 1, let
 $\cup_{j=1}^{k(1)}\{ a_j^1 \}$ be a finite subset
of $\partial B_2\setminus\f_0 (X)$
fulfilling $(2_1)$ (for the set $A$ in Lemma 1 take $\partial B_2\setminus
\f_0 (X)$).
Consider the compact set $K:=\bar B_1 \cup \f_0 (X_{R_0})$.
By Lemma~\ref{klem:2}
in \cite{FGR} the polynomially convex hull $\hat K$ of $K$ is contained in
$\bar B_1 \cup \f_0 (X)$; in particular it does not contain any of the points
$a_j^1$.
 We can choose distinct points $x^1_1,x^1_2,\ldots,x^1_{k(1)}$
in $X\setminus ( \hat K \cup X_{R_0})$ such that
$$\rho (x^1_i)>\max_{x\in \f_0 ^{-1} (B_1)} \rho (x)
\quad i=1,2,\ldots, k(1)$$  and use Proposition~\ref{prop:f1}
 $k(1)$ times to find an
automorphism $\Phi^1_1$ with $\Vert \Phi^1_1 (z)-z\Vert\le \frac{\epsilon_1}{
2}$
for $z\in K$, $\Phi^1_1 \circ \f_0 (x_j^1) = a_j^1$ $j=1,2,\ldots, k(1)$ and
$\Phi^1_1(z)=z+O(\vert z-b_i\vert^{k+1}$ as $z\to b_i$ $i=1,2,\ldots,n+2$.
By lemma~\ref{klem:2} we find another automorphism $\Phi^2_1$ which moves the
compact set
$\Phi^1_1(K)$
less than $\frac{\epsilon_1}{ 2}$, fixes the points $a_j^1$ $j=1,2,\ldots,
k(1)$,
matches the identity up to order $k$ at the points $b_i$ $i=1,2,\ldots,n+2$
and has the property that the submanifold $\Phi^2_1\circ \Phi^1_1\circ \f_0 (
\tilde X)$ does not osculate of order $k$ at any point
$\Phi^2_1\circ\Phi^1_1\circ\f_0 (x)$ with
$x \in (X_{ R_{0}}\cup \tilde X)\setminus U_1$.
The composition $\Phi_1:=\Phi_1^2\circ\Phi_1^1$ (together with the set
$\cup_{i=1}^{k(1)} \{x_i^1 \}$ satisfies all properties from
$(3_1)$ to $(8_1)$. Finally  we choose $R_1$ big enough to satisfy $(9_1)$.

The description of the $m$-th step is similar to the first step. To be
accurate we carry it out in detail. Suppose we have already constructed
$\epsilon_i, R_i>0$,
 finite subsets $\cup_{j=1}^{k(i)} \{a_j^i\}$
of $\partial B_{i+1}$, finite subsets $\cup_{j=1}^{k(i)}\{x_j^i\}$ of $X$
together with  automorphisms $\Phi_i$ satisfying $(1_i)$ up to $(9_i)$
for all $i$ from $1$ to $m-1$.
Again we first choose $\epsilon_m\le \epsilon_{m-1}$ so small that any
perturbation of the embedding $\Phi_{m-1}\circ\f_0 :\tilde X \ha \C^n$
which is smaller than $3\epsilon_m$ on the compact set
$(X_{ R_{m-2}}\cup \tilde X)\setminus U_{m-1}$ does not destroy the property
that $\Phi_{m-1}\circ\f_0 (\tilde X)$ does not osculate of order $k$ at any
point in the image of $(X_{ R_{m-2}}\cup \tilde X)\setminus U_{m-1}$.
Next, using Lemma 1, we choose a finite
subset $\cup_{j=1}^{k(m)} \{a_j^m \}$ of $\partial B_{m+1}\setminus
\Phi_{m-1}\circ f_0 (X)$ fulfilling $(2_m)$.
Consider the compact set $K:=\bar B_m \cup \Phi_{m-1}\circ\f_0 (X_{R_{m-1}})$.
Again by lemma~\ref{klem:2}
in \cite{FGR} the polynomially convex hull $\hat K$ of $K$ is contained in
$\bar B_m \cup \Phi_{m-1}\circ\f_0 (X)$. Hence it
does not contain any of the
points $a_j^m$.
We  choose distinct points $x^m_1,x^m_2,\ldots,x^m_{k(m)}$
in $X\setminus ( \hat K \cup X_{R_{m-1}})$ such that
$$\rho (x^m_i)>\max_{x\in \Phi_{m-1}\circ\f_0 ^{-1} (B_m)} \rho (x)
\quad i=1,2,\ldots, k(m)$$  and use  Proposition~\ref{prop:f1}
 $k(m)$ times to find an
automorphism $\Psi^1_m$ with $\Vert \Psi^1_m (z)-z\Vert\le \frac{\epsilon_m}{
2}$
for $z\in K$, $\Psi^1_m\circ \f_0 (x_j^m) = a_j^m$ $j=1,2,\ldots, k(m)$,
$\Psi^1_m (a_j^l) = a_j^l$ $j=1,2,\ldots, k(l)$ $\forall \ l\le m-1$,
 and $\Psi^1_m
(z)=z+O(\vert z-b_i\vert^{k+1})$ as $z\to b_i$ $i=1,2,\ldots,n+2$.
By lemma~\ref{klem:2}
we find another automorphism $\Psi^2_m$ which moves the compact set
$\Psi^1_m (K)$
less than $\frac{\epsilon_m}{ 2}$, fixes the points $a_j^l$ $j=1,2,\ldots,
k(l)$
$1\le l\le m$,
matches the identity up to order $k$
at the points $b_i$, $i=1,2,\ldots,n+2$,
and has the property that the submanifold
$\Psi^2_m\circ \Psi^1_m\circ \Phi_{m-1}\circ \f_0 (\tilde X)$ does not osculate
of order $k$ at any point in the image of
$ (X_{ R_{m-1}}\cup \tilde X)\setminus U_m$.
We set $\Psi_m:=\Psi_m^2\circ\Psi^1_m$ and $\Phi_m=\Psi_m\circ \Phi_{m-1}$.
Now all conditions
$(3_m)$ to $(8_m)$ are also satisfied. Finally we choose $R_m$ big enough to
satisfy $(9_m)$ and $(10_m)$.

By Proposition 2,
 the properties $(5_m)$ together with the fact that $\epsilon_m
< \frac{1}{ 3^m}$ imply that $\lim_{m\to \infty} \Phi_m$ exists (uniformly on
compacts) on $\Omega:=\cup_{m=1}^\infty \Phi_m^{-1} (B_m) $ and defines a
biholomorpic map from $\Omega$ onto $\C^n$.

If $z\in\f_0 (X)$, then we find an $m$ such that
$x:=\f_0^{-1} \in X_{R_m}$. So properties $(4_m)$ imply:
$$\Vert \Phi_l \circ \f_0 (x) - \Phi_m\circ\f_0 (x)\Vert \le
\sum_{k=m}^{l-1} \Vert \Phi_{k+1}\circ\f_0 (x) -\Phi_k\circ\f_0 (x)\Vert
\le \sum_{k=m}^{l-1} \epsilon_{k+1} <\sum_1^\infty \epsilon_i <\infty$$
This shows that the set $\{\Phi_m (z)=\Phi_m\circ\f_0 (x): m\in \N\}$ is
bounded in $\C^n$. By Proposition 2 this means $z\in \Omega$. We have proved
 $\f_0 (X) \subset \Omega$. This shows that
$\f:=\Phi\circ \f_0$ is an (proper holomorphic) embedding
$\f: X\ha \C^n$.

We now show that $\f$ satisfies the conclusion of the theorem.
First observe that $(6_m)$ for all $m$ implies
$$\f (x_j^m) = a_j^m \quad j=1,2,\ldots,k(m) \quad\forall \ m\in\N. $$
We set $\alpha=\sum_{i=1}^\infty \epsilon_i$. By $(1_m)$ we have
$\alpha<\frac{1}{ 2}$.
Next we show:
$$\f ^{-1} (B_{m-1})\subset (\Phi_m\circ\f_0)^{-1} (B_{m}) \quad \forall \
m\in\N\hfill\eqno{(\star)}$$

Let $x\in \f^{-1} (B_{m-1})$ be an arbitrary point. Since this means
$\Phi \circ \f_0 (x) \in B_{m-1}$,
we may choose $k_0>m$ and $\epsilon=1-2\alpha >0$,
such that
$$\Phi_k \circ \f_0 (x) \in B_{m-1+\epsilon} \quad\forall \ k\ge k_0
\hfill\eqno{(\star\star )}.$$
By $(5_{m+1})$ and Rouche's theorem (see \cite{Ch} p.110) we have
$\Psi_{m+1} (B_m)\supset B_{m-2\epsilon_m}$, i.e.,
$B_m\supset \Psi_{m+1}^{-1} (B_{m-2\epsilon_m})$. Hence
$\Phi_m^{-1} (B_m) \supset \Phi_m^{-1} (\Psi_{m+1}^{-1} (B_{m-2\epsilon_m}))=
\Phi_{m+1}^{-1} (B_{m-2\epsilon_m})$. By induction, using $(5_{m+2}),\ldots ,
(5_k)$ we find
$$
\Phi_m^{-1} (B_m)\supset \Phi_k^{-1} (B_{m-2\sum_{l=m}^{k-1}
\epsilon_l})\supset \Phi_k^{-1} (B_{m-2\alpha}).
$$
Together with
$(\star\star )$ this
implies by our choice of $\epsilon$ that $\f_0 (x) \in \Phi_m^{-1}
(B_m)$
and hence
$x\in (\Phi_m\circ\f_0)^{-1} (B_{m}) $. So $(\star)$ is proved.

Now suppose $F:\C^n\to\C^n$ is a nondegenerate holomorphic map with
 $F^{-1}(\C^n\setminus \varphi (X))= \C^n\setminus \varphi(X)$ and
 $\varphi^{-1}\circ F\circ\varphi:X\to X $ is the element $\alpha_t:=\alpha
(t,\cdot)$
of the given family of selfmaps of $X$.
By moving the origin by an arbitrarily small translation,
 we can assume $JF (0)\ne 0$.
We set $\beta =\prod_{i=1}^\infty (1-\epsilon_i)>0$
and choose $m_0$ big enough so that for all
$m\ge m_0$ it follows that
$$
t\in T_m,\quad JF (0)> \frac{1}{ m^n \beta^n},\quad
F(0)\in B_{\frac{m+3}{2}} \hfill \eqno{(A)}
$$
We fix a number $m\ge m_0$ and claim that $$\forall \ i\ge m+2\quad:
 F(B_m)\cap\cup_{j=1}^{k(i)}\{ a_j^i\} =\emptyset\hfill \eqno{(B)}
$$

Suppose the contrary, i.e., there exists $z\in B_m$ with $F(z)=a_j^i$ for some
$i\ge m $ and $1\le j\le k(i)$. Since $F(\C^n\setminus\f (X))\subset
\C^n\setminus \f (X)$, we have $z\in \f (X)$.
Let $x=\f^{-1} (z)\in \f^{-1} (\B_m)$, i.e.,
$F\circ \f (x)=a_j^i=\f (x_j^i)$. This means $\alpha_t (x) = \f^{-1}\circ
F\circ \f (x)= x_j^i $.
But by $(\star)$ we have $x\in (\Phi_{m+1}\circ\f_0)^{-1} (B_{m+1})$
and,
 since
$t\in T_m $, it follows that
$$\rho (\alpha_t (x))\le
\max_{y\in (\Phi_{m+1}\circ \f_0)^{-1} (\bar \B_{m+2}), t \in T_{m+2}}
\rho (\alpha (t,y)).$$
Hence, since $i-1\le m+1$, we have
$$\rho (\alpha_t (x))\le
\max_{y\in (\Phi_{i-1}\circ \f_0)^{-1} (\bar \B_{i}), t\in T_{i}}
\rho (\alpha (t,y)) .$$
According to $(3_i)$ this, together with the choice of $x_j^i$, implies
$\rho (x_j^i) > \rho (\alpha_t (x) ) $ which contradicts $x_j^i=\alpha_t (x)$.
Thus property $(B)$ is proved.

There exists a natural number $k$ such that
$F (B_m) \subset B_{k+2}$.  Suppose $k>m+2$. We consider the holomorphic maps
$$F_j :\B_1 \to B_{k+2},\quad F_j(z):= F(z \cdot m \prod_{l=j+1}^k
(1-\frac{\epsilon_l}{ 2})), \quad j=m+2,\ldots,k$$
For $j=m+2,\ldots, k$ we have:

First, $F_j (\B_1) =F( B_{m \prod_{l=j+1}^k (1-\frac{\epsilon_l}{ 2})})
\subset F(B_m)$. Hence, from $(B)$ it follows that
 $F_j (B_1)$ does not contain any point
$a_l^j$, $l=1,2,\ldots,k(j)$.

Second, $JF_j (0) =m^n\
\prod_{l=j+1}^k (1-\frac{\epsilon_l}{ 2})^n JF (0) > m^n \beta^n
JF(0)>1$
is a consequence of
$(A)$.

And third, also by $(A)$, we have $F_j (0)=F(0)\in B_\frac{j+1}{ 2}$.

Using $(1_k)$  for the holomorphic map $F_k$,
these three properties imply,
$$
F_k (B_{(1-\frac{\epsilon_k}{ 2})})=F_{k-1} (B_1)\subset B_{k+1}.
$$
Proceeding
by induction from $k$ down to $m+2$
and using property  $(1_j)$
for the holomorphic map $F_j$
we find in the same way
that $ F_j (B_1) \subset B_{j+2}$.
 For $j=m+2$ this means $$F (B_{m\cdot \beta}) \subset
F_{m+2} (B_1)\subset B_{m+4}\quad \forall\ m\ge m_0.$$
This growth restriction forces $F$ to be an affine map.
Since $F$ is a non-degenerate map, it must be an affine
automorphism.

From $(9_m)$ and the fact that $\f_0 (\tilde X)$ osculates of order $k$ at
the points $\f_0 (x_i)=b_i$, $i=1,2,\ldots,n+2$ it follows that the
submanifold $\f (\tilde X)$ osculates of order $k$ at the points $\f (x_i)$,
$i=1,2,\ldots,n+2$. From $(8_m)$ (and the accurate choice of $\epsilon_m$ in
the beginning of the m-th step) it follows that $\f (\tilde X)$ does
not osculate
of order $k$ at any other point.
Since the affine automorphism $F$ leaves
$\f (X)$ invariant, it is clear that $\f (\tilde X)$ (the union  of
the components of maximal dimension of the smooth part of $\f(X)$)
is $F$-invariant as well.
Furthermore an affine
automorphism preserves the property of
osculating of order $k$ at some point. Thus
 $F$ permutes the  points $b_1,b_2,\ldots,b_{n+2}\in \C^n$. But
no affine automorphism except the identity permutes these points. Hence
$F=\id$.
\end{proof}


\begin{thebibliography}{Bla}
\newcommand{\bibauthor}{\em}
\newcommand{\bibjour}{\relax}

\bibitem{AM}
{\bibauthor S.S. Abhyankar, T.T. Moh},
Embeddings of the line in the plane.
{\bibjour  J.~Reine u. Angew.~Math. }\bf 276\rm, (1975), 148--166

\bibitem{A}
{\bibauthor D. Akhiezer},
Lie group actions in complex analysis.
Aspects of Mathematics E27.
Viehweg 1995.

\bibitem{AL}
{\bibauthor E. Andersen, L. Lempert},
  On the group of
 holomorphic automorphism of $\C^n$
{\bibjour  Invent. Math. }\bf 110 \rm  (1992), 371--388

\bibitem{BF}
{\bibauthor G.T. Buzzard, J.E. Forn\ae ss},
An embedding of $\C$ into $\C^2$ with hyperbolic complement.
{\bibjour  Math. Ann. }\bf 306 \rm  (1996), 539--546

\bibitem{BFo}
{\bibauthor G.T. Buzzard, F. Forstneric},
A Carleman type theorem for proper
holomorphic embeddings
{\bibjour  Ark. Mat }\bf 35 \rm (1997), 157--169

\bibitem{Ch}
{\bibauthor E. Chirka},
Complex analytic sets.
Kluwer, Dordrecht 1989

\bibitem{F1}
{\bibauthor  F. Forstneric},
Interpolation by holomorphic automorphisms and embeddings in $\C^n$.
Preprint.(1996),

\bibitem{F}  
{\bibauthor F. Forstneric},
 Actions of $(\R,+)$ and $(\C,+)$ on complex manifolds.
{\bibjour  Math.~Z. }\bf 223 \rm(1996), no. 1, 123--153

\bibitem{FGR}
{\bibauthor F. Forstneric, J. Globevnik, J.P. Rosay},
Nonstraightenable complex lines in $ \C^2$.
{\bibjour  Ark. Mat. }\bf 34 \rm (1996), no. 1, 97--101.

\bibitem{FR}
{\bibauthor F. Forstneric, J.P. Rosay},
Approximation of biholomorphic mappings
by automorphisms of $\C^n$.
{\bibjour  Invent. Math. }\bf 112 \rm (1993), 323--349

\bibitem{GK}
{\bibauthor P. Griffiths, J. King},
Nevanlinna Theory and holomorphic Mappings between Algebraic
Varieties.
{\bibjour  Acta Math }\bf 130\rm, (1973), 145--220


\bibitem{H}
{\bibauthor W.Y. Hsiang},
Cohomology Theory of Topological Transformation Groups.
Springer 1975

\bibitem{I} {\bibauthor S. Iitaka},
On logarithmic Kodaira dimension of algebraic varieties.
{\sl in}
Complex Analysis and Algebraic Geometry.
175--189. Iwanami Shoten, Tokyo 1977.


\bibitem{J} {\bibauthor Z. Jelonek},
The extension of regular and rational embeddings.
{\bibjour  Math.~Ann. }\bf 277\rm, (1987), 113-120

\bibitem{K}
{\bibauthor S. Kaliman},
Extensions of isomorphisms between affine algebraic subvarieties of
$k^n$
to automorphisms of $k^n$.
{\bibjour  Proc.~A.M.S. }\bf 113\rm, (1991), 325--334

\bibitem{RR1}
{\bibauthor J.P. Rosay, W. Rudin},
Holomorphic maps from $\C^n$ to
$\C^n$. {\bibjour  Trans. A.M.S. }\bf 310\rm, (1988), 47--86

\bibitem{RR2}
{\bibauthor J.P. Rosay, W. Rudin},
Holomorphic embeddings of $\C$ in $\C^n$. 563--569
{\sl in} Several Complex Variables.
Math.~Notes 38.
Proceedings of the Mittag Leffler Institute 1987-88.
Princeton University Press 1993

\bibitem{Sa}
{\bibauthor F. Sakai},
Kodaira dimension of complements of divisors.
{\sl in}
Complex Analysis and Algebraic Geometry.
239-257. Iwanami Shoten, Tokyo 1977.

\bibitem{Sch}
{\bibauthor J. Sch\"urmann},
Embeddings of Stein spaces into affine spaces of minimal
dimension.
{\bibjour  Math. Ann. }\bf 307\rm, no.3 , (1997), 381--399

\bibitem{S}
{\bibauthor V. Srinivas},
On the embedding dimension of the affine variety.
{\bibjour  Math. Ann. }\bf 289\rm, (1991), 125--132

\bibitem{JW1}
{\bibauthor J. Winkelmann},
On Free Holomorphic $\C $-actions on $\C^n$ and
Stein Homogeneous Manifolds.
{\bibjour  Math. Ann. }\bf 286\rm, (1990), 593--612

\bibitem{JW2}
{\bibauthor J. Winkelmann},
On Automorphisms of Complements of Analytic Subsets in $\C^n$.
{\bibjour  Math. Z. }\bf 204\rm, (1990), 117--127

\end{thebibliography}
\end{document}